%% file: Paper_GP_arxiv.tex
\documentclass{article}
\usepackage{arxiv}

\usepackage{color} 
\usepackage{epstopdf}%
\usepackage{subfigure}%
\usepackage{verbatim}
\usepackage{tikz}
\tikzset{every picture/.style={/utils/exec={\sffamily}}}
\usetikzlibrary{mindmap,trees}
\usepackage{pgfplots}
\usepgfplotslibrary{fillbetween}
\pgfplotsset{compat=newest}
\usepackage{multirow}
\usepackage{mathcomp}
\usepackage{rotating}
\usepackage{booktabs}
\usepackage{textcomp}
\usepackage{amsmath}
\usepackage{amssymb}
\usepackage{svg}
\usepackage{amsbsy}

\usepackage[utf8]{inputenc}
\usepackage[T1]{fontenc}
\usepackage{graphicx}
\usepackage{dsfont}
\usepackage{amsmath}
\usepackage{mathtools}
\usepackage{amssymb}
\usepackage[prependcaption]{todonotes}
\usepackage[noend,ruled]{algorithm2e}
\SetKwIF{If}{ElseIf}{Else}{if}{}{else if}{else}{end if}%
\SetKwInput{KwData}{Input}
\SetKwInput{KwResult}{Output}
\usepackage{nicefrac}
\usepackage{tabularx}
\usepackage{multirow}
\usepackage{siunitx} %
\usepackage{eurosym}
\usepackage{color}
\usepackage{pgfplots}
\usepackage{booktabs}
\usepackage{subfigure}
\usepackage{verbatim}
\usepackage{csquotes}
\usetikzlibrary{patterns}
\usepgfplotslibrary{fillbetween}
\usepackage{rotating}
\usepackage{mathcomp}

\usepackage{standalone}
\usepackage{longtable}

\usepackage[prependcaption]{todonotes}

\newcommand{\objective}[1]{&\min \span &\rlap{$\displaystyle #1$}}

\newcommand{\basicmass}{M_{\text{0}}}
\newcommand{\batterymass}{M_{\text{b}}}
\newcommand{\passengermass}{M_{\text{d}}}
\newcommand{\additionalmass}{M_{\text{a}}}
\newcommand{\totalmass}{m}

\newcommand{\cycle}{\Lambda}
\newcommand{\speed}{v}
\newcommand{\slope}{s}
\newcommand{\timet}{t}

\newcommand{\prob}{\pi}

\newcommand{\motortorque}{t^{\text{M}}}
\newcommand{\wheeltorque}{T^{\text{W}}}
\newcommand{\transmissionratio}{i}
\newcommand{\transmissionratioOne}{i_1}
\newcommand{\transmissionratioTwo}{i_2}
\newcommand{\motorrotspeed}{\omega^{\text{M}}}
\newcommand{\motortorquemax}{\overline{t}^{\text{M}}}

\newcommand{\motorpowerfactor}{k^{\text{M,p}}}
\newcommand{\motorpowerlossone}{p^{\text{M,L,1}}}
\newcommand{\motorpowerlosstwo}{p^{\text{M,L,2}}}
\newcommand{\motorpowerlossthree}{p^{\text{M,L,3}}}
\newcommand{\motorpowerlossconst}{p^{\text{M,L,const}}}
\newcommand{\motorpowerout}{p^{\text{M}}}
\newcommand{\motorspeedn}{n^{\text{M}}}
\newcommand{\motorspeednmax}{\overline{N}^{\text{M}}}
\newcommand{\motorpowerin}{p^{\text{M,in}}}
\newcommand{\motorpowerintotal}{{p}^{\text{M,in,avg}}}
\newcommand{\motorpoweroutmax}{\overline{p}^{\text{M,out}}}
\newcommand{\motormass}{m^{\text{M}}}

\newcommand{\bin}[2]{b_{#1,#2}}

\newcommand{\wheelspeedn}{N^{\text{W}}}
\newcommand{\transmissionratiomin}{\underline{I}}
\newcommand{\transmissionratiomax}{\overline{I}}
\newcommand{\motorpowerlossconstant}{P_{\text{L,ref,const}}}

\newcommand{\R}{\mathbb{R}^+}

\DeclareRobustCommand{\stirling}{\genfrac\{\}{0pt}{}}

\begin{document}

\title{Efficient Powertrain Design --- A Mixed-Integer Geometric Programming Approach}

\author{Philipp Leise \\ 
Technische Universit\"{a}t Darmstadt \\
Chair of Fluid Systems \\
Department of Mechanical Engineering \\
\texttt{philipp.leise@fst.tu-darmstadt.de}\\
\\
 \And 
 Peter F. Pelz \\
 Technische Universit\"{a}t Darmstadt \\
Chair of Fluid Systems \\
Department of Mechanical Engineering \\
\texttt{peter.pelz@fst.tu-darmstadt.de}\\}
\maketitle

\begin{abstract}
The powertrain of battery electric vehicles can be optimized to maximize the travel 
distance for a given amount of stored energy in the traction battery. 
To achieve this, a combined control and design problem has to be solved which 
results in a non-convex Mixed-Integer Nonlinear Program. 
To solve this design task more efficiently, we present a new systematic optimization approach 
that leads to a convex Mixed-Integer Nonlinear Program. 
The solution process is based on a combination of Geometric Programming and a Benders decomposition.
The benefits of this 
approach  are a fast solution time, a global convergence, and the ability to 
derive local sensitivities in the optimal design point with no extra cost, as 
they are computed in the optimization procedure by solving a dual problem.
The presented approach is suitable for the evaluation of a complete driving cycle, 
as this is commonly done in powertrain system design, or for usage in a stochastic 
approach, where multiple scenarios are sampled from a given probability density 
function. The latter is useful, to account for the uncertainty in the driving 
behavior to generate solutions that are optimal in an average sense for a 
high variety of vehicles and drive conditions.  
For the powertrain model we use a transmission model with up to two selectable 
transmission ratios 
and an electric motor model that is based on a scaled efficiency map representation. 
Furthermore, the shown model can also be used to 
model a continuously variable transmission to show beneficial energy savings 
based on this additional degree of freedom.
The presented design approach is also applicable to a wide variety of design tasks for technical 
systems besides the shown use case.
\end{abstract}

\section{Introduction}

To improve the air quality, especially in cities, and to enable a more sustainable 
transportation, battery electric vehicles (BEVs) are a promising technology.  
With rapidly decreasing investment costs for the traction battery, as shown by 
\cite{nykvist2015rapidly}, the profitability of these vehicles increases and therefore 
the number of new registrations. 
As a result the market share of electric vehicles increased over the last years, as 
shown by \cite{palmer2018total}. This increase of new registrations of BEV is also 
driven by environmental protection regulations in major markets, like for example 
China, \cite{li2016consumers}. 
The major drawback of BEVs in comparison to vehicles with an internal combustion 
engine, as stated by \cite{egbue2012barriers}, is a shorter travel distance.
The lower ratio of stored energy to weight within the vehicles and the higher recharging 
times of BEVs in comparison to fuel powered vehicles lead to requirements to improve 
the overall efficiency of the powertrain as much as possible. 

The optimization program for the design task to derive optimized powertrain configurations must also consider 
the usage phase, cf.~\cite{Silvas.2016b}, to derive meaningful system designs.
Therefore, a combined component design and control problem has to be solved.
Furthermore, we explicitly consider a multi-speed transmission to improve the system efficiency even further.
This problem results in a non-convex Mixed-Integer Nonlinear Program (MINLP) which is hard to solve 
to global optimality. 
As a result, often primal heuristics like a genetic algorithm (GA) or particle-swarm optimization are 
used to derive a (local-optimal) optimized solution. For this approach see e.g. \cite{schonknecht2016electric},
\cite{tan2018gear} and \cite{walker2013modelling}.

To derive global-optimal solutions more rapidly, we present a new approach that combines \emph{Geometric Programming} (GP) and 
a \emph{Benders Decomposition} to solve the given design task. This results in a separation of the overall non-linear and non-convex 
optimization program in a main Mixed-Integer Linear Program and a convex Nonlinear Program (NLP) as a subproblem.
These are then solved in an iterative solution process, where the convex subproblem leads in each iteration 
to additional constraints 
within the main problem. The main problem, on the other hand, leads to new binary assignments within the subproblem. 
This iterative solution of smaller optimization problems in combination with the convexified system properties 
results in an increased solution speed in comparison to the original non-convex MINLP.  

In the following, we present a brief overview about geometric programming and convex 
powertrain design. Afterwards, we present the new geometric programming model 
for a powertrain design, the developed Benders approach  and a heuristic to derive good solutions fast.
 Then, we show an outlook on 
the usage of mixed-integer geometric programming in powertrain design and close with a conclusion 
and summary.

\section{Geometric Programming}

The presented optimization model relies on the geometric programming approach, 
which was first developed by Duffin, Peterson, and Zener in 1967, cf.~\cite{duffin1967geometric}. 
It is a modeling approach in which a $\log$-transformation is used to convert a 
non-convex nonlinear program in a convex program. 
This transformation is possible, as shown in \cite{boyd2007tutorial}, if the 
objective and constraints either consist of monomials
\begin{equation}
m(x) = \alpha_0 x_1^{\alpha_1} x_2^{\alpha_2} \dots x_n^{\alpha_n},
\end{equation}
or posynomials 
\begin{equation}
p(x) = \sum_k \alpha_{0k} x_{1k}^{\alpha_{1k}} x_{2k}^{\alpha_{2k}} \dots x_{nk}^{\alpha_{nk}}.
\end{equation}
In this example the $n$ variables are given by $x$, while $\alpha$ represents specific parameters 
for each constraint and $k$ the total number of terms of the sum. 
The underlying optimization program can be transformed in a convex program by using a $\log$-transformation 
if it has a monomial as an objective and monomial or posynomial constraints of the following form, \cite{boyd2007tutorial} :
\begin{equation}
\begin{aligned}
\min_{x}& \quad & m(x) \\
\textrm{s.t.}& \quad & m(x) &= 1  \\
 & &p(x) &\leq 1  \\
 & &x &> 0. \\
\end{aligned}
\end{equation}
The resulting $\log$-transformed problem is convex. In general, convex programs can be solved 
very efficiently by state-of-the-art solvers. We refer to \cite{boyd2007tutorial} as a general 
introduction on geometric programming and to \cite{boyd2004convex} for an overview about 
convex optimization.
The geometric programming approach is currently used in aircraft design, as for example 
shown by \cite{hoburg2014geometric, Burton.2017}.

Besides the geometric programming approach shown below, further research was done to 
derive convex optimization models for the powertrain component sizing and control. For instance 
\cite{murgovski2012component} present a detailed modeling approach that relies on a convex 
optimization program. 
Additional approaches for convex modeling and optimization are given in \cite{murgovski2013engine}, 
\cite{egardt2014electromobility} and \cite{hu2016greener}.
To the best of the author's knowledge, the powertrain design and control problem for 
battery electric vehicles with an explicit transmission ratio design has not yet been modeled 
as a (mixed-integer) geometric program. 
Hence, we present in the following section an easily extendable modeling approach, 
which is also compliant with the underlying physical principles.

\section{Preprocessing}

The powertrain of battery electric vehicles primarily consists of the battery, a power 
electronics, an electric motor and a transmission with one or more transmission ratios. For simplicity, 
we first present only the integration of a transmission with one variable single transmission 
ratio and an example electric motor. Further components are currently neglected, but can 
be added later, if necessary.
In an early design stage in powertrain development commonly a backwards dynamic 
longitudinal vehicle model is used to derive the torque and speed requirements at 
the wheel, \cite{tran2020thorough}. 
In the following we rely on a given driving load representation $\cycle$:
\begin{equation}
\cycle = \begin{bmatrix}
\mathbf{\speed}^T\\
\mathbf{\slope}^T  \\
\boldsymbol{\prob}^T
\end{bmatrix}, \cycle \in \mathbb{R}^{3\times l}.
\label{eq:dataset}
\end{equation}
Here, $\mathbf{v}$ represents the vector of velocity requirements and $\mathbf{s}$ 
represents the slope requirements.
This representation is equivalent to a given driving cycle representation, cf.~\cite{Silvas.2016b}, 
if the additional weights $\prob_t$, $t~\in~(1,\dots, l)$ are neglected.  
In the following we represent the set of used discrete time steps with $\mathcal{T}$.
We can also interpret the given data representation as a finite event set approximation 
of a stochastic program with recourse. 
In this case an uncertain variable is given by the representation of a drive conditions. 
This interpretation enables us to evaluate multiple vehicle representations and drive 
conditions at once by using for instance specific weights $\prob_t$ and not a driving cycle 
but a set of representative scenarios that approximate the underlying probability density 
function of the uncertain drive conditions. 
This set of scenarios is usually smaller than the complete cycle.
This approach is useful in many situations, as most time steps in the driving cycle do not 
affect the final optimal solution.
The usage of a comparable stochastic approach is for instance shown by \cite{caillard2014optimization} 
or \cite{wasserburger2020risk}.

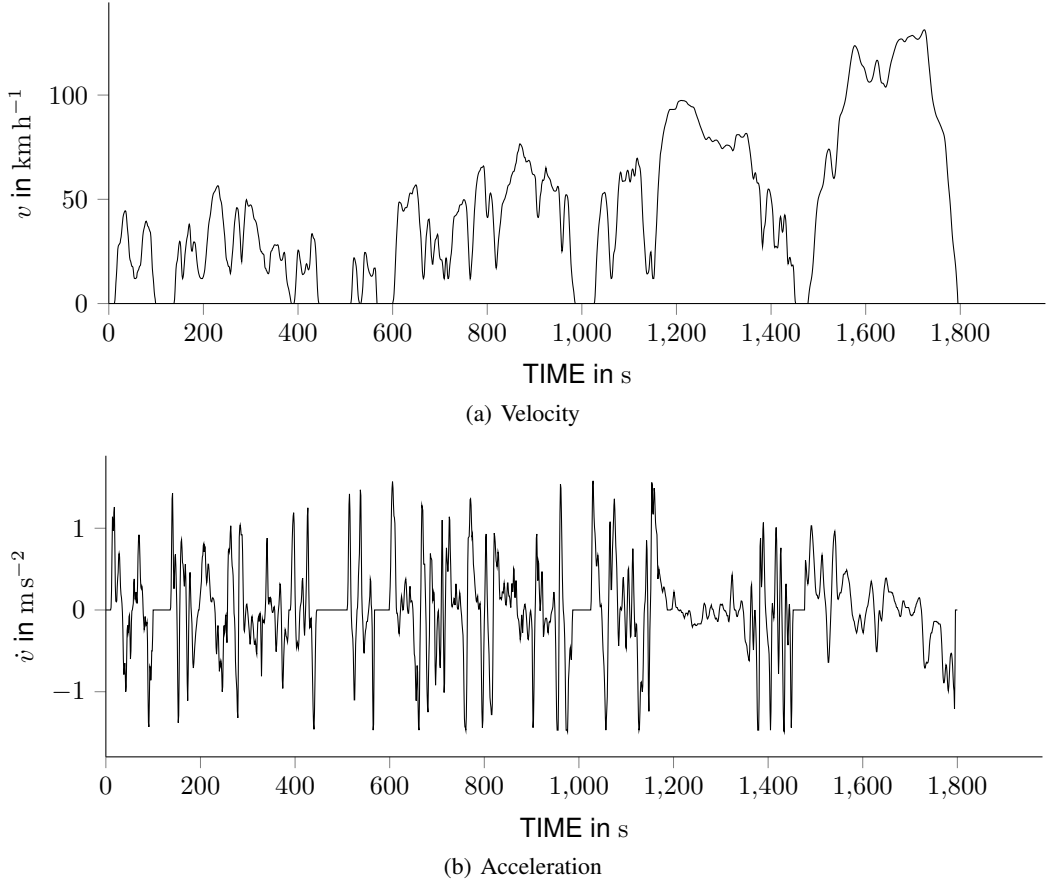
\begin{figure}
\center
\subfigure[Velocity]{
\begin{tikzpicture}
\begin{axis}
[xlabel near ticks,
  xlabel={TIME in \SI{}{\second}},
  ylabel near ticks,
  ylabel={$v$ in \SI{}{\kilo \meter \per \hour}},
  xtick align=outside,
  ytick align=outside,
  axis x line*=bottom,
  axis y line*=left,
  xmin=0,
  ymin=0,
  legend pos= north east,
  legend style={draw=none},
  legend cell align=left,
  height = 0.2\paperheight,
  width = 0.65\paperwidth
  ]
  
\addplot[line legend] 
	table[x=t,y=v,col sep=comma]{WLTC.csv};
\end{axis}
\end{tikzpicture}
}
\subfigure[Acceleration]{
\begin{tikzpicture}
\begin{axis}
[xlabel near ticks,
  xlabel={TIME in \SI{}{\second}},
  ylabel near ticks,
  ylabel={$\dot{v}$ in \SI{}{\meter \per \square \second }},
  xtick align=outside,
  ytick align=outside,
  axis x line*=bottom,
  axis y line*=left,
  xmin=0,
  legend pos= north east,
  legend style={draw=none},
  legend cell align=left,
  height = 0.2\paperheight,
  width = 0.65\paperwidth
  ]
  
\addplot[line legend] 
	table[x=t,y=a,col sep=comma]{WLTC.csv};
\end{axis}
\end{tikzpicture}
}
\caption{Input based on WLTC, \cite{tutuianu2013development}. (a) velocity requirements. (b) acceleration 
requirements derived by using a central difference scheme.}
  \label{fig:wltc}
\end{figure}

We use the following longitudinal vehicle model for modeling the torque and speed 
requirements in the cycle 
\cite[p. 83ff.]{mitschke2014dynamik}:
\begin{equation}
 T_t = \left( \lambda_{\text{i}}\,\totalmass\dot{v}_t +  \totalmass g\left(\sin\left(\alpha\right) 
 + \lambda_{\text{r}} \right) + \frac{1}{2}\rho c_{\text{w}} A v_t^2 \right)r \quad \forall t\, \in\, \mathcal{T}.
\label{eq:torquegearbox}
\end{equation}
We refer to Tab.~\ref{tab:parameterVehicle} for the meaning and values of the 
used parameters.
\begin{table}
\center
\caption{Used vehicle parameters for dimensioning the powertrain.}
\begin{scriptsize}
\begin{tabular}{clll}
\toprule
 PARAMETER & VALUE & UNIT & DESCRIPTION \\
\midrule
$\basicmass$ & $1100$    & \SI{}{\kilo \gram} & vehicle mass\\
$\additionalmass$ & $75$  &  \SI{}{\kilo \gram} & additional component mass in powertrain\\
$\passengermass$ & $75$  &  \SI{}{\kilo \gram} & mass of the considered driver\\
$\batterymass$ & $550$ & \SI{}{\kilo \gram} & mass of the used battery\\
$g$ & $9.81$ & \SI{}{\meter \per \second \squared} & specific gravitational constant \\
$\alpha$ & $0$ & -- & terrain slope angle \\
$\rho$ & $1.2041$ & \SI{}{\kilo \gram \per \cubic \meter} & air density at \SI{20}{\celsius} and sea level\\
$c_{\text{w}}$ & 0.3 & -- & drag coefficient \\
$A$ & 2.2 & \SI{}{\meter \squared} & vehicle reference area \\
$\lambda_{\text{r}}$ & $0.01$ & -- & rolling resistance coefficient \\
$r$ & 0.3 & \SI{}{\meter} & wheel radius \\
$\lambda_i$ & 1.0 & -- & inertia factor \\
$\eta_{g}$ & $0.98$ & -- & efficiency of the gearbox\\
$\motorspeednmax$ & $10000$ & \SI{}{\per \minute} & maximum speed of the selected motor\\
\end{tabular}
\end{scriptsize}
\label{tab:parameterVehicle}
\end{table}
For the longitudinal vehicle model, we only use the speed $\speed$ and acceleration $\dot{\speed}$ in 
each time step in the given driving cycle as inputs. 
We derive acceleration values for each given discrete time step $\timet \in \mathcal{T}$ in the cycle 
based on a central difference scheme. For simplicity, we assume $\slope_t = 0 \, \forall \, \timet \in \mathcal{T}$. 
Nevertheless, this assumption can be easily dropped, as the modeling approach is independent of the 
slope values. It only differs in the preprocessing.
As a reference cycle, we use in this section the \emph{worldwide harmonized
light vehicles test cycle}  (WLTC), \cite{tutuianu2013development}, as shown in Fig.~\ref{fig:wltc}. 
Furthermore, we only consider positive acceleration ($\dot{v} > \SI{0}{\meter\per\square\second}$) 
and speed ($v > \SI{0}{\meter\per\second}$) values in the following. 
We also present a further model extension to integrate recuperation in section~\ref{sec:extensions}.

As the power demand in the cycle depends on the total mass $\totalmass$ of the vehicle, 
we can not estimate the power directly before optimization, if the used mass of relevant 
components are modeled as variables within the model. 
For a first approach, we use in this example a fixed value which is estimated as 
follows: $\totalmass~=~\basicmass + \passengermass + \additionalmass + \batterymass$. 
Here, $\basicmass$ is the basic vehicle mass, $\passengermass$ the mass for the considered 
passengers, $\additionalmass$ the additional mass of powertrain components, and $\batterymass$ 
the battery mass.
Furthermore, we show a possible extension to a variable mass in section \ref{sec:extensions}.

\begin{figure}[ht]
	\centering
  \includegraphics[
   width=0.95\textwidth,
   trim={0cm 0cm 0cm  0cm}, 
   clip, 
   scale=1,
  ]{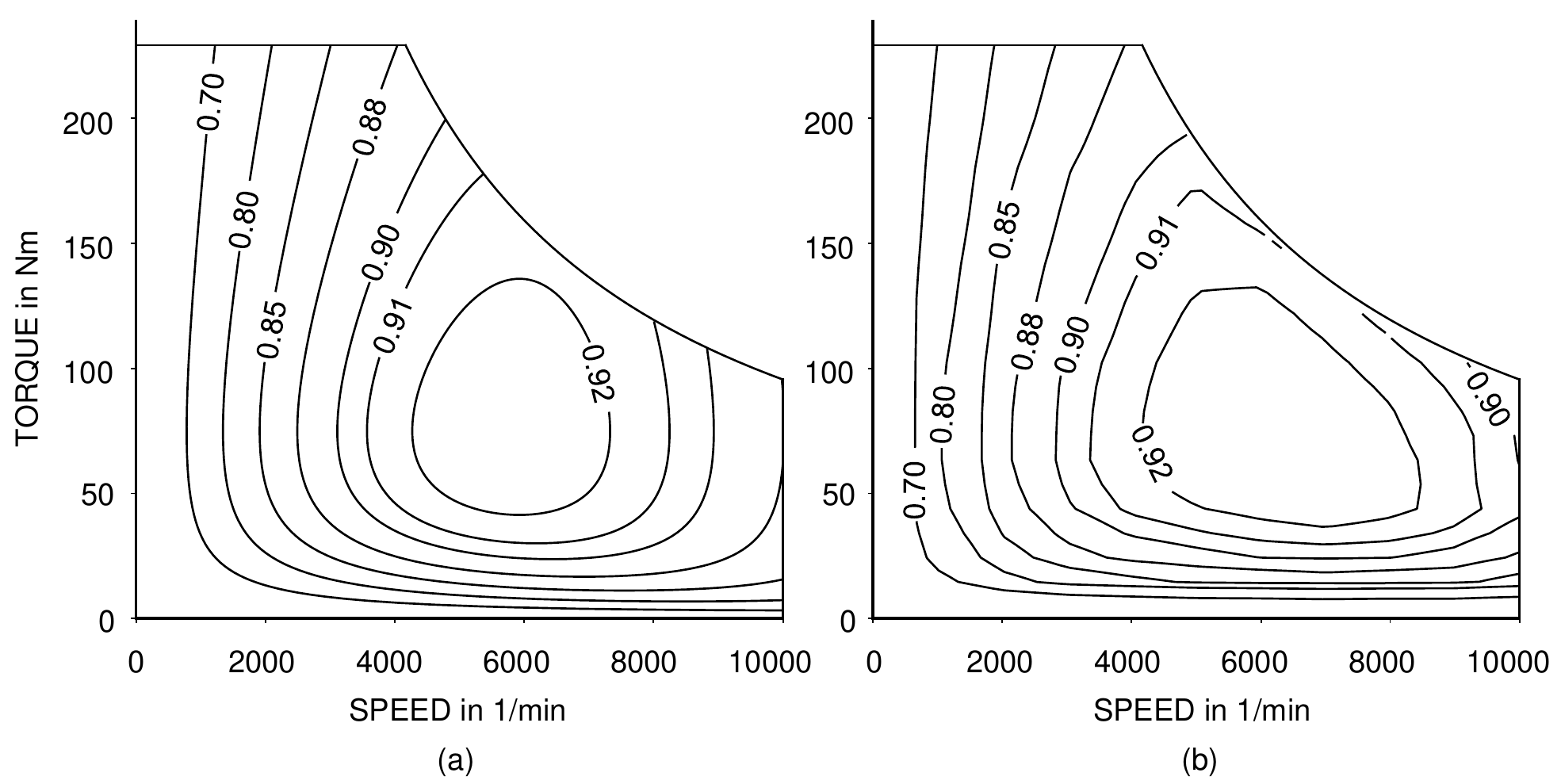}
	\caption{Comparison of used model and reference for the first quadrant of the 
  efficiency map. The given values represent isolines with the same efficiency. 
  (a) Output of the used GP-compatible efficiency map model for a motor with \SI{100}{\kilo\watt}. 
  (b) Reference efficiency map of a given scaled motor with the same power, cf.~\cite{An.2017}. }
	\label{fig:motorefficiencymap}
\end{figure}

Within BEV, permanent magnet synchronous motors (PMSM) are commonly used, since they have a weight 
and efficiency benefit compared to other motor technologies. Therefore, the underlying electric 
motor model used in this contribution represents a PMSM. We use multiple monomials to represent 
specific power losses of the motor. For instance \cite{vratny2019conceptual, mahmoudi2015efficiency} 
showed that PMSM partial losses, like the friction loss in bearings, or stray losses can be 
approximated by using monomials. 
As a reference for comparison, we use a scaled version of an efficiency map shown by \cite{An.2017} 
for a motor developed for the usage in electric vehicles.
For simplicity, we currently do not model the power electronics and battery explicitly. 
Instead, we only use an aggregated efficiency map and assume that it represents the motor as 
well as the required power electronic components.
We scale the derived efficiency map within the optimization according to the motor's maximum torque $\motortorquemax$.
This scaling is useful in the early design stage and for instance also shown in \cite{vratny2019conceptual}.
It assumes that the efficiency map of a derived electric motor is comparable with the 
efficiency map of a reference motor of the same product series. 
The benefit of the shown approximation method is that it is compliant with the GP approach, 
as it is modeled with the help of multiple monomial and posynomial constraints.
An example result of the used motor efficiency approximation method is shown in 
Fig.~\ref{fig:motorefficiencymap} for the \SI{100}{\kilo\watt} reference motor.

Next to the given driving cycle this optimization method allows to integrate easily further 
constraints, as for instance high speed requirements or starting constraints based on the gradeability on a slope. 
They can also be interpreted as further scenarios / time steps that must be met but have 
a zero probability $\pi_t$. 
This ensures that they are not accounted for in the objective but have to be fulfilled by the final solution. 

\section{Used Basic Model}

The complete model is shown in description \ref{prob:GP} in a condensed form.
The used parameters, sets and variables are given in Tab.~\ref{tab:variables}. 
\begin{table}%
\center
  \caption{GP model notation}
  \label{tab:variables}
  \begin{scriptsize}
    \begin{tabularx}{\textwidth}{@{}clll}
    \toprule
    SET &  &  & DESCRIPTION \\
\midrule
$\mathcal{T}$ & & & set of all considered time steps \\
$\mathcal{T}^+$ & & & set of all steps based on the given cycle\\
\midrule
PARAMETER & VALUE & UNIT & DESCRIPTION \\
\midrule
$\eta_g$ & $0.98$ & -- & transmission efficiency of gearbox\\  
$A$ & $0.4161$  & -- & motor hyperbola constant\\
$\underline{I}$ & $1.0$  & -- & minimum transmission ratio\\
$\overline{I}$ & $18.0$  & -- & maximum transmission ratio\\
$\pi$ & $1/\sum_{k=1}^l 1$ & -- & probability of occurrence for all $|\mathcal{T}^+|$ time steps in the given cycle \\
$\wheeltorque$ & based on cycle & \SI{}{\newton \meter} & wheel torque \\
$\wheelspeedn$ & based on cycle & \SI{}{\per \minute} & wheel rotational speed\\
$\motorspeednmax$ & $10000$ & \SI{}{\per \minute} & motor maximum rotational speed \\
$P_{\text{ref}}$ & 100 & \SI{}{\kilo\watt} & power of reference motor\\ 
$P_{\text{L, ref, 1}}$ & $787.35$ & \SI{}{\watt} & power reference loss 1 \\ 
$P_{\text{L, ref, 2}}$ & $1566.67$ & \SI{}{\watt} & power reference loss 2 \\ 
$P_{\text{L, ref, 3}}$ & $9904.85$& \SI{}{\watt} & power reference loss 3\\ 
$P_{\text{L, ref, const}}$ & $1059.34$ & \SI{}{\watt} & power reference constant loss \\ 
$\overline{P}^{\text{M, in}}$ & $8000$ & \SI{}{\kilo\watt} & upper limit of required motor power \\
\\
      VARIABLE & DOMAIN & UNIT & DESCRIPTION \\
      \midrule
$\transmissionratio$ & $[\transmissionratiomin, \transmissionratiomax]$  & -- & transmission ratio \\  
$\motortorque$ & $[0, \motortorquemax]$ & \SI{}{\newton \meter} & motor torque \\
$\motortorquemax$ & $\R$ & \SI{}{\newton \meter} & motor maximum torque \\
$\motorspeedn$ & $[0,\motorspeednmax]$ & \SI{}{\per \minute} & motor speed \\
$\motorpowerintotal$ & $\R$ & \SI{}{\watt} & average used power in cycle \\
$\motorpowerin$ & $\R$ & \SI{}{\watt} & motor power requirement drawn from battery \\
$\motorpoweroutmax$ & $\R$ & \SI{}{\watt}  & motor maximum design power \\
$\motorpowerout$ & $\R$ & \SI{}{\watt}  & motor design power \\
$\motorpowerfactor$ & $\R$ & --  & motor power factor
\end{tabularx}
      
\end{scriptsize}
\end{table}
\begin{figure}[htp]
\label{prob:GP}
	\begin{subequations}
	\begin{align}
	\objective{\motorpowerintotal}\label{eq:obj}\\
	&\;\mbox{s.t.}&& \sum_{t\in\mathcal{T}^+}  \prob_t \motorpowerin_t \leq \motorpowerintotal && \label{eq:mod2}\\
	&&& \transmissionratiomin \leq \transmissionratio \leq \transmissionratiomax && \label{eq:mod2a}\\
	&&&  \motortorque_t \transmissionratio\, \eta_g = \wheeltorque_t &&  \forall t \in \,\mathcal{T}\label{eq:mod3}\\
	&&&  \motorspeedn_t = \transmissionratio \wheelspeedn_t &&  \forall t \in \,\mathcal{T}\label{eq:mod4}\\
		&&&  \motorpoweroutmax =\motorspeednmax \frac{2\pi}{60} A \motortorquemax&& \label{eq:mod8}\\
	&&&  \motorspeedn_t \frac{2\pi}{60} \motortorque_t \leq \motorpoweroutmax  &&  \forall t \in \,\mathcal{T}\label{eq:mod5}\\
	&&& \motorspeedn_t \leq \motorspeednmax &&  \forall t \in \,\mathcal{T}\label{eq:mod6}\\
	&&& \motortorque_t \leq \motortorquemax &&  \forall t \in \,\mathcal{T}\label{eq:mod7}\\
	&&&  \motorpowerfactor =\frac{\motorpoweroutmax}{P_{\text{ref}}} &&  \label{eq:mod9}\\
	&&&  \motorpowerlossone_t = \motorpowerfactor P_{\text{L, ref, 1}} \left(\frac{\motortorque_t}{\motortorquemax}\right)\left(\frac{\motorspeedn_t}{\motorspeednmax A}\right)^{3.93}  &&  \forall t \in \,\mathcal{T}\label{eq:mod10}\\
	&&&  \motorpowerlosstwo_t = \motorpowerfactor P_{\text{L, ref, 2}} \left(\frac{\motortorque_t}{\motortorquemax}\right)  &&  \forall t \in \,\mathcal{T}\label{eq:mod11}\\
		&&&  \motorpowerlossthree_t = \motorpowerfactor P_{\text{L, ref, 3}} \left(\frac{\motortorque_t}{\motortorquemax}\right)^{2}  &&  \forall t \in \,\mathcal{T}\label{eq:mod12}\\	
		&&&  \motorpowerlossconst_t = \motorpowerfactor \motorpowerlossconstant  &&  \forall t \in \,\mathcal{T}\label{eq:mod13}\\	
		&&& \motorpowerout_t = \motortorque_t \frac{2\pi}{60} \motorspeedn_t	&&  \forall t \in \,\mathcal{T} \label{eq:mod14}\\
		&&& \motorpowerin_t \geq \motorpowerout_t + \motorpowerlossone_t + \motorpowerlosstwo_t +  \motorpowerlossthree_t + \motorpowerlossconst_t &&  \forall t \in \,\mathcal{T} \label{eq:mod15} \\
		&&& \motorpowerin_t \leq \overline{P}^{\text{M, in}} &&  \forall t \in \,\mathcal{T} \label{eq:mod16}
		\end{align}
\end{subequations}
\caption{GP model}
\end{figure}
The model reads as follows: As an objective, \eqref{eq:obj}, we consider an upper bound of the average 
approximation of the total power losses in all time steps $\mathcal{T}^+$, which is given in constraint \eqref{eq:mod2}. 
The transmission ratio $\transmissionratio$ must lay between an upper bound $\transmissionratiomax$ 
and a given lower bound $\transmissionratiomin$, \eqref{eq:mod2a}. 
Constraints \eqref{eq:mod3} and \eqref{eq:mod4} model the transmission with the variable transmission 
ratio $\transmissionratio$ and the fixed gearbox efficiency $\eta_g$. We use the speed in \SI{}{\per \minute} 
of the given motor. If desired, it is also possible to use the angular speed $\motorrotspeed = 2\pi \motorspeedn/{60}$ instead.
As the motor can be scaled in the torque direction, we model the variable motor power with \eqref{eq:mod8} 
and the upper hyperbola limit of the feasible motor domain with constraint \eqref{eq:mod5}.
All torque values $\motortorque$ and speed values $\motorspeedn$ must lay in the feasible motor 
domain. Next to the previously modeled hyperbola, we also model the upper torque and speed 
limits in \eqref{eq:mod6} and \eqref{eq:mod7}.
As the efficiency map of the underlying PMSM should be scaled according to the maximum torque, 
we use a scaling parameter $\motorpowerfactor$ and introduce it in a monomial in constraint \eqref{eq:mod9}. 
We use the motor power $ P_{\text{ref}}$ as a reference. Constraints \eqref{eq:mod10} -- \eqref{eq:mod13} represent specific power losses, based on the derived efficiency map model. We used a nonlinear least-square method to estimate the reference values and exponents.
We use the scaling factor $\motorpowerfactor$  in each loss term to get an according 
scaling of all losses.
The constant losses in constraint \eqref{eq:mod13} for instance can represent air gap 
losses and cooling power requirements, cf.~\cite{vratny2019conceptual}.
These losses could be estimated with further scaling laws as well. Nevertheless, we 
omit these further details, as they don't affect the shown optimization approach and 
could be added later, if desired. 
The shaft power in each time step is given in \eqref{eq:mod14} and the total required 
power that must be supplied by the battery is given in a relaxed form in \eqref{eq:mod15}. 
We additionally model an upper bound for $\motorpowerin$ in \eqref{eq:mod16}. 
This upper bound is not active in the final result, but required for solving the GP.
The complete model represents the physical properties of the given system.
But especially the posynomial constraints (\ref{eq:mod2}, \ref{eq:mod15}) have to be 
relaxed to allow the usage in the GP. Nevertheless, in the final solution they meet the 
equality condition to derive a system design that has the lowest power losses. Therefore, 
these constraint relaxations do not affect the found solution.

\section{Exemplary Results for the Basic Model}

We use the weights $ \pi_t = 1/\sum_{k=1}^n 1 \, \forall t \in \mathcal{T}^+$ for each step 
based on the set $\mathcal{T}^+$ which was derived from the driving cycle. Here, we only 
consider positive speed and acceleration values, which results for the WLTC in $767$ time steps. 
We additionally use two further restrictions, which must be fulfilled by the optimized powertrtain.
First a start torque at speed $\motorspeedn = \SI{0.2}{\per \minute}$ that represents a 
gradeability on a slope of \SI{66}{\percent}. The speed value is set arbitrarily above zero to enable a GP solution. 
As the  selected low speed value lays below the motor design speed, it does not affect 
the solution. 
Furthermore, we require a high-speed limit of \SI{160}{\kilo\meter\per\hour} with an 
acceleration of  \SI{0.3}{\meter\per\square\second}. 
These two constraints are modeled with two additional time steps, that must be fulfilled 
in the final motor domain, but are not used in the objective. 
In total, these $769$ time steps result in approx. $6100$ variables and approx. $20800$ 
constraints for the complete program.
The final transmission which is based on the WLTC, the required additional loads, and the given 
vehicle parameters in Tab.~\ref{tab:parameterVehicle} has a transmission ratio of $7.06$. 
The selected motor has a maximum power $\motorpoweroutmax$ of \SI{185}{\kilo\watt} which 
results in a maximum torque of $\motortorquemax \approx \SI{424}{\newton\meter}$. 
These values are comparable with a solution derived by a genetic algorithm with a higher 
fidelity model that uses the original motor efficiency map next to recuperation.

Besides the specific transmission ratio $\transmissionratio$ in the single-speed transmission, we can also 
consider a continuously variable transmission (CVT). In this case the transmission ratio is not 
modeled as a single variable but instead as a vectorized variable for all time steps. This 
additional degree of freedom results in a more efficient system design and can be seen as a 
physical lower bound for multi-speed transmission. We used an upper bound $\transmissionratiomax = 18$ 
besides the lower bound $\transmissionratiomin = 1.0$. 
In this domain, the optimized system design uses a motor with \SI{73}{\kilo\watt} and a 
maximum torque $\motortorquemax \approx \SI{167}{\newton\meter}$. 
The CVT result also reduces the objective value compared to the single-speed transmission 
by approx. \SI{8}{\percent}. This value depends on the given upper limit $\transmissionratiomax$. 
Nevertheless, we can conclude that a CVT can improve the overall system performance. 
The usage of multiple transmission ratios in a BEV powertrain, which can be interpreted as a discrete 
approximation of the CVT, can result in potential energy savings, cf.~\cite{schonknecht2016electric}.

We showed that the developed GP model is suitable for the early design stage of powertrain design.
With the given model it is possible to compute optimized powertrain design solutions rapidly 
in comparison to a non-convex modeling or commonly used heuristics.
We developed the model and preprocessing in Python. Therefore, we used GPkit \cite{burnell2020gpkit} 
for modeling and Mosek 9.2 as a solver. The computations were done on a Debian-based 
machine with an Intel i5-5200U and 12 GB RAM. On this hardware, we were able to compute 
the optimized system designs in less than three seconds each. This emphasizes the high 
potential for usage in rapid development in an early design stage.

\section{Mixed-Integer Geometric Programming}

Besides the already shown CVT and single-speed transmission, we can also derive solutions for multi-speed transmissions 
with this new design methodology. 
Therefore, we extend the shown CVT model in \eqref{prob:GP}, by using two additional variables for the considered
transmission ratios. Here, $\transmissionratioOne$ represents the first possible transmission ratio, while 
$\transmissionratioTwo$ represents the second possible transmission ratio. 
Then, we can think of the final optimal solution having 
an optimal binary vectorized mapping $\bin{t}{g}$ for each time step $t \in \mathcal{T}$ and 
each transmission ratio $g \in \mathcal{G}$ with $\mathcal{G} =  \{1,2\}$.
This binary decisions then indicate if either the first or second transmission ratio is chosen within the 
considered time step.
This approach can also be extended to transmissions with more than two transmission ratios. For simplicity, we only restrict our 
description on the case of two transmission ratios. 

The resulting MINLP has two stages. In the first stage, we assign the aforementioned continuous design variables.
In the second stage, we then assign the aforementioned binary control variables for each time step to achieve a global optimal solution.
As a result, we get a control strategy for the vehicle, based on the underlying dataset. 

The second stage decision to assign transmission ratios to each time step is a combinatorial 
subproblem, that depends on the results of the first stage. 
This coupling between both stages in combination with the non-convex property within the original 
design space are the major challenges within the solution process of the underlying MINLP.

In the following, we want to focus on the underlying combinatorial problem in the second stage and consider its 
complexity properties to derive approaches to improve the solution speed.
As the assignment of a single transmission ratio to each time step (single-speed solution) is 
trivial, we will only consider transmissions 
that have at least $2$ transmission ratios ($|\mathcal{G}| \geq 2$). 
We always assume, that $|\mathcal{G}| < |\mathcal{T}|$, as the assignment is otherwise (CVT solution) again trivial.
All possibilities of an assignment of transmission ratios $g \in \mathcal{G}$ to time steps are 
described by the cartesian product, which results in 
\begin{equation}
N_{g \rightarrow t} = |\mathcal{G}|^{|\mathcal{T}|}.
\label{eq:cartprod}
\end{equation}
Nevertheless, this fast growing number of possible assignments, $N_{g \rightarrow t}$, can be 
reduced considerably, if we consider explicitly symmetry influences. %

Symmetries in Integer Programs and Mixed-Integer Programs represent permutable variable assignments in the 
final solution, which lead to equivalent objective values and solutions, that are equivalent in terms 
of the final solutions properties. These symmetries lead to major difficulties for current state-of-the-art solvers 
that use a branch-and-cut approach as for example stated by \cite{margot2010symmetry}. These symmetries lead to 
a high number of equivalent enumerations within the solution process, which can reduce the solution speed considerably.

In our case, the presented MINLP has a high number of symmetric solutions, as the assignment of transmission ratios $g \in \mathcal{G}$ 
to each time step $t \in \mathcal{T}$ is equivalent, if we permute the assigned ratios in all time steps accordingly. 
The final solutions for each permuted assignment is equivalent, as each transmission ratio is optimized individually from each other. 
The final number of combinatoric solutions, after removing all symmetries, can be represented, 
in our case, by stirling numbers of the second kind, cf.~\cite[p.~243ff.]{graham1990concrete} and \cite[p.~66ff.]{knuth1997art}. 
These are used in combinatorics to derive the number of possibilities to assign a set of $n$ distinguishable 
objects, in our case the time steps, into $k$ nonempty indistinguishable subsets, in our case the transmission ratios. 
Stirling numbers of second kind, $\stirling{n}{k}$, are defined by:
\begin{equation*}
  S_{n,k} = \stirling{n}{k} \coloneqq \frac{1}{k!}\sum\limits_{i=0}^{k}(-1)^{i}\binom{k}{i} \left(k - i\right)^n.
\end{equation*}
The final number of distinguishable combinatoric possibilities, without redundant symmetric assignments, 
is then considerably smaller, than proposed by Eq. \eqref{eq:cartprod}.
For the considered showcase of two-speed transmission the Stirling number then reduces to 
\begin{equation}
S_{n,2} = \stirling{n}{2} = 2^{n-1} -1.
\end{equation}

\subsection{Brute Force Reference} 

One possibility for a small problem size is the usage of a brute force approach to solve all combinatoric 
possibilities. A fixed transmission ratio to time step assignment then reduces the MINLP to an NLP which
can then be solved again by using the aforementioned approach. 

To increase the solution speed of the brute force approach, we only compute assignments that are not symmetric, as stated before.
This can be achieved with a lexicographic ordering of all binary assignments before optimization of each instance.
We only consider the first half of all assignment vectors in the lexicographic ordering, 
as the second half corresponds to symmetric problem definitions. 
It is important to mention that this is only valid for two-speed transmissions. 
For three- or more-speed transmission the subset with all removed symmetric binary 
assignments has to be chosen differently.

Since even the Stirling number increases rapidly with the number of time steps, a 
brute force approach is only suitable for small problem instances with 
a low number of time steps. 
To recall, this small number of time steps can also be interpreted as specific 
scenarios, if the weights $\pi_t$ within Eq.~\eqref{eq:dataset} are used within optimization.
This reduced approach then is also often able to derive comparable results to the usage of a complete cycle.

Within this article, we only focus on the derived new design and solution methodology. 
We use the derived brute force solutions to evaluate the applicability of the 
Benders decomposition and heuristic solution approach more systematically.

\subsection{Benders Decomposition}

As a second solution approach, besides the brute force approach, we present the usage 
of a Benders decomposition, cf.~\cite{benders1962partitioning} and \cite{rahmaniani2017benders}. 
It is suitable for the problem at hand, since it allows the exploitation of the specific problem structure given 
by the afore mentioned reference model. 
Within this approach, the MINLP is 
divided into two problems. One main problem in which all binary control decisions are selected and a subproblem in which 
all continuous design decisions are set. These problems are then solved iteratively and within 
each iteration specific information between both problems is shared. 
Here, the main problem is further constraint by so-called \emph{optimality cuts}, which are derived from dual information 
within the subproblem. After reoptimization of the main problem a new binary assignment vector is transferred 
to the subproblem. After a new optimization this leads to new dual variable assignments in the intermediate solution and then to 
further new optimality cuts within the main problem.     

\subsection{Heuristic Solution Approach}

As a third approach, we also present a heuristic, which uses an iterative assignment strategy. 
It uses the physical property that it is more efficient for a transmission ratio --- time step assignment that 
low speed and high torque loads are assigned to the first transmission ratio, while high speed and 
low torque loads are assigned to the second transmission ratio. With this domain-specific knowledge, it is 
then possible to assign the time steps iteratively to each transmission ratio. 
The solution procedure starts with the CVT solution and then iteratively adds additional constraints for loads, 
starting with the highest torque and the highest speed loads. Since a fixation of the transmission ratios for specific loads leads usually to 
an infeasibility of the former CVT solution, we reoptimize the new model with its additional constraints.
Then, the remaining unassigned loads in each time step are iteratively assigned based on the chosen transmission ratios 
and the difference to the chosen already fixed transmission ratios. 

\begin{figure}
	\centering
  \includegraphics[
   width=0.95\textwidth,
   trim={0cm 0cm 0cm  0cm}, 
   clip, 
   scale=1,
  ]{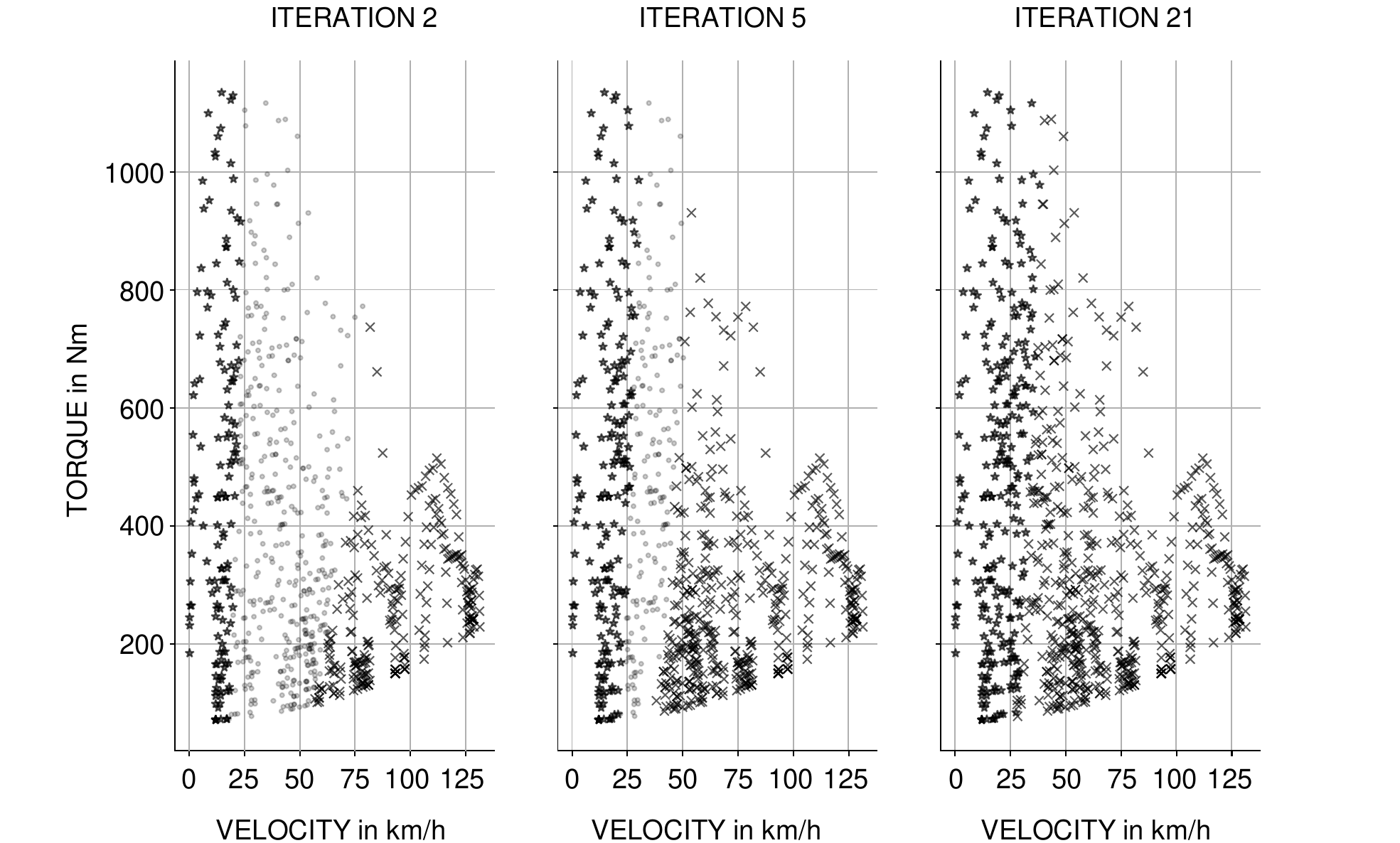}
	\caption{Result of the developed heuristic to derive 
  near-global-optimal solutions for the WLTC. A star marks the usage of the first transmission ratio, 
  while a cross shows the second transmission ratio. Unset binaries are marked with a circle. 
  After $21$ iterations all binary decision variables are set in this example.}
	\label{fig:heursistic}
\end{figure}

The result of the developed heuristic for the WLTC example, which we already used as a previous example, 
with 2 transmission ratios is shown in Fig.~\ref{fig:heursistic}.
It is able to derive an assignment for all loads within $21$ iterations. The most assignments are done within the 
first iterations, while later iterations lead to the assignment of only a few remaining unset binary decisions.

\section{Comparison of the Developed Solution Approaches}

In the following, we present first results to show the applicability of the selected approaches. 
We compare the results and solution procedures of the Benders decomposition approach, the brute force approach, and the 
heuristic approach. Therefore, we used $6$ -- $12$ scenarios in conjunction with the already shown model (\eqref{prob:GP}).
To assure a more general understanding of the algorithms, we used $120$ different instances.
Here, we used the legislative driving cycles WLTC, Artemis Motorway, FTP75 and NYCC. Additionally, 
we use one cycle which is based on real driving behavior. It was derived by members of the 
institute of mechatronics systems at Technische Universität Darmstadt and we name 
this cycle ``IMS Alltracks''. We refer to \cite{Esser.2018} for more details on the synthesis process. 
The scenarios were derived by using a k-means clustering approach based on the given datasets. 

Within this section, we only focus on the algorithmic benefits and do not discuss any technical results that can be derived from the 
usage of the developed optimization approaches. This ensures a more concise overview on the algorithmic capabilities.

An overview of all instances is given in the Appendix in Table~\ref{tab:a1}. 
The Benders decomposition approach is able to derive in all  $120$ test-cases the global-optimal solution, as resulting from the brute force approach.
Additionally, a combination of the Benders decomposition and the heuristic is also able to derive the same global optimal solution.
Nevertheless, the solution time can not be significantly improved by combining both, 
heuristic and Benders decomposition, on the given test-set. 
The Benders approach finds all solutions in a fraction of the required iterations of the brute force approach. 
Therefore, it is suitable for a more rapid solution procedure that can also derive global optimal solutions.
The developed heuristic is also able to often derive the global optimal solution. Nevertheless, 
sometimes the resulting solution is only near-global optimal. But this near global optimality can also be effective for practical 
usage, where the solution time is a significant restriction.   

Hence, we can conclude that the shown approaches are effective to derive global (in case of Benders decomposition)
 or at least near-global (in case of the heuristic) optimal solutions.
We will evaluate the abilities on larger datasets and with an extended level of detail of the underlying model later. 
Therefore, we present further model extensions that can be used within the NLP (GP) and MINLP (MIGP) approach.

\section{Model Extensions}
\label{sec:extensions}

With the help of the shown GP model we were able to transform the non-convex original 
problem to a convex program, without the need of further polynomial or piecewise-affine 
model approximations of components or the preliminary selection of transmission ratios.
The underlying motor representation is based on physical properties of the modeled PMSM.
The shown model \eqref{prob:GP} can be extended in multiple ways.
For a more detailed mass estimation besides the shown fixed mass, we can model the 
total mass $\totalmass$, based on variable partial masses of different components.
For instance, we can model the motor mass $\motormass$ with a specific mass $\rho^{\text{M}}$ 
and the maximum motor power $\motorpoweroutmax$, as shown in constraint~\eqref{eq:mod17}. 
This results in a new variable mass $\totalmass$ which can be estimated with fixed 
reference masses for the basic vehicle $\basicmass$, for additional powertrain 
components $\additionalmass$, for the passengers $\passengermass$,
and for the battery $\batterymass$ \eqref{eq:mod18}. 
Furthermore, the torque at the wheel now depends on the mass, which is linear 
in $\totalmass$ and can be approximated in a posynomial by using \eqref{eq:mod19}.
\begin{subequations}
\begin{align}
		&&& \motormass = \rho^{\text{M}} \motorpoweroutmax \label{eq:mod17}\\
		&&& \totalmass \geq \basicmass + \passengermass + \additionalmass + \batterymass + \motormass \label{eq:mod18}\\
		&&& t_t^{\text{W}} \geq \delta_t \totalmass  + \gamma_t &&  \forall t \in \,\mathcal{T} \label{eq:mod19}
\end{align}
\end{subequations}

Additionally, we can easily extend the given model to account for recuperation. 
For this step, we have to model the fourth quadrant of the torque-speed domain 
with a second efficiency map model. 
This second model is comparable to the given model of the first quadrant, 
but differs slightly in the given parametrization. We can add the recuperation 
domain by using the absolute values of the given torque to map it to a new 
positive domain. In this domain, we use the newly defined efficiency map model 
and integrate it in the final model with a second weighted sum in 
constraint~\eqref{eq:mod2}. To account for the relation between the driving and 
recuperation, we add further positive weights, which weigh both sides against each other.

The shown GP-approach is also suitable for the already mentioned stochastic representation. 
In this, we can use multiple scenarios and weights to represent the given 
uncertain driving conditions. 
To ensure robust final solutions of the optimized system that fulfill all 
desired drive conditions in the underlying representative driving cycle or 
probability density function, we have to add further loads that lay on the 
convex hull of each cluster or at least at the highest required power demands. 
They are added as constraints but do not affect the objective.

Within this contribution, we  presented a white-box reference-based model for the motor efficiency.
Next to this approach further more detailed component models can be cast as a geometric program.
For instance a more detailed motor mass estimation besides the density based approach would be possible.
Furthermore, each component model representation which can be represented as a 
quantity of monomial or posynomial constraints can be used. This is often possible by adding further variables as shown in the motor model.
If it is not possible, we can also fit a posynomial approximation, as for instance shown by \cite{hoburg2014fitting}.

\section{Summary and Conclusion}

We presented a mixed-integer geometric programming model for the early design stage of battery 
electric vehicles. We showed that the optimal design of a given gearbox and electric 
motor can be represented as a convex program, without specific requirements on the 
variables and constraints. This easily extendable basic model can therefore be 
used in numerous possibilities.
The presented model was developed based on publicly available software and the 
usage of an of-the-shelf optimizer. Therefore, it is easily applicable in industry and academia. 
System model details next to the motor efficiency and transmission ratios can be 
added easily to the shown model, which makes it very suitable for early design 
proof-of-concepts and the comparison of optimal powertrain designs with 
different system layouts or components.
Additionally, we showed an extension to explicitly consider  a 
discrete transmission ratio
selection and an according control strategy. 
Within the future, we will extend this model with an explicit consideration 
of recuperation and 
evaluate the shown method in more detail with more driving cycles, larger datasets and potentially other 
optimization methods.
Furthermore, this optimization approach can also be transferred to other domains of designing technical systems.
Therefore, we present an example for the design of a water supply system that is based on \cite{leise2018energy}.

%\paragraph*{Author contributions}
%P.L.: methodology, implementation, writing; P.P.: supervision

\bibliographystyle{unsrt}  
\bibliography{biblio}

\newpage
\section*{Appendix}

\input{table_comparison.tex}

\end{document}

%% file: table_comparison.tex
\renewcommand{\thetable}{A1}

\begin{scriptsize}
\begin {longtable}{ccccccccc}%
    \caption{We present in the following all instances ordered by its identification number (ID) that were used for the comparison of the Benders (B), heuristic (H) and brute force (BF)
    approach to solve the developed Mixed-Integer Geometric Program. B and BF lead to a gap of $0$. 
    All primal solutions (P) are given in the logarithmic GP space. The number of used iterations is given by ``\#it''. } \\
ID&BF/P& B/P& B+H/P & H/P & B/\#it& B+H/\#it & BF/\#it & CYCLE\\%
\hline
\pgfutilensuremath {1}&\pgfutilensuremath {9.87}&\pgfutilensuremath {9.87}&\pgfutilensuremath {9.87}&\pgfutilensuremath {9.87}&\pgfutilensuremath {20}&\pgfutilensuremath {16}&\pgfutilensuremath {127}&Artemis Motorway\\%
\pgfutilensuremath {2}&\pgfutilensuremath {9.78}&\pgfutilensuremath {9.78}&\pgfutilensuremath {9.78}&\pgfutilensuremath {9.78}&\pgfutilensuremath {27}&\pgfutilensuremath {33}&\pgfutilensuremath {127}&Artemis Motorway\\%
\pgfutilensuremath {3}&\pgfutilensuremath {9.95}&\pgfutilensuremath {9.95}&\pgfutilensuremath {9.95}&\pgfutilensuremath {9.95}&\pgfutilensuremath {23}&\pgfutilensuremath {25}&\pgfutilensuremath {127}&Artemis Motorway\\%
\pgfutilensuremath {4}&\pgfutilensuremath {9.72}&\pgfutilensuremath {9.72}&\pgfutilensuremath {9.72}&\pgfutilensuremath {9.72}&\pgfutilensuremath {16}&\pgfutilensuremath {22}&\pgfutilensuremath {127}&Artemis Motorway\\%
\pgfutilensuremath {5}&\pgfutilensuremath {9.89}&\pgfutilensuremath {9.89}&\pgfutilensuremath {9.89}&\pgfutilensuremath {9.89}&\pgfutilensuremath {25}&\pgfutilensuremath {26}&\pgfutilensuremath {127}&Artemis Motorway\\%
\pgfutilensuremath {6}&\pgfutilensuremath {9.79}&\pgfutilensuremath {9.79}&\pgfutilensuremath {9.79}&\pgfutilensuremath {9.79}&\pgfutilensuremath {31}&\pgfutilensuremath {36}&\pgfutilensuremath {127}&Artemis Motorway\\%
\pgfutilensuremath {7}&\pgfutilensuremath {9.3}&\pgfutilensuremath {9.3}&\pgfutilensuremath {9.3}&\pgfutilensuremath {9.3}&\pgfutilensuremath {20}&\pgfutilensuremath {22}&\pgfutilensuremath {127}&FTP75\\%
\pgfutilensuremath {8}&\pgfutilensuremath {8.6}&\pgfutilensuremath {8.6}&\pgfutilensuremath {8.6}&\pgfutilensuremath {8.6}&\pgfutilensuremath {31}&\pgfutilensuremath {32}&\pgfutilensuremath {127}&FTP75\\%
\pgfutilensuremath {9}&\pgfutilensuremath {9.42}&\pgfutilensuremath {9.42}&\pgfutilensuremath {9.42}&\pgfutilensuremath {9.42}&\pgfutilensuremath {19}&\pgfutilensuremath {21}&\pgfutilensuremath {127}&FTP75\\%
\pgfutilensuremath {10}&\pgfutilensuremath {9.44}&\pgfutilensuremath {9.44}&\pgfutilensuremath {9.44}&\pgfutilensuremath {9.44}&\pgfutilensuremath {29}&\pgfutilensuremath {25}&\pgfutilensuremath {127}&FTP75\\%
\pgfutilensuremath {11}&\pgfutilensuremath {9.17}&\pgfutilensuremath {9.17}&\pgfutilensuremath {9.17}&\pgfutilensuremath {9.17}&\pgfutilensuremath {27}&\pgfutilensuremath {26}&\pgfutilensuremath {127}&FTP75\\%
\pgfutilensuremath {12}&\pgfutilensuremath {8.92}&\pgfutilensuremath {8.92}&\pgfutilensuremath {8.92}&\pgfutilensuremath {8.92}&\pgfutilensuremath {30}&\pgfutilensuremath {32}&\pgfutilensuremath {127}&FTP75\\%
\pgfutilensuremath {13}&\pgfutilensuremath {9.12}&\pgfutilensuremath {9.12}&\pgfutilensuremath {9.12}&\pgfutilensuremath {9.14}&\pgfutilensuremath {18}&\pgfutilensuremath {17}&\pgfutilensuremath {127}&IMS Alltracks\\%
\pgfutilensuremath {14}&\pgfutilensuremath {9.24}&\pgfutilensuremath {9.24}&\pgfutilensuremath {9.24}&\pgfutilensuremath {9.24}&\pgfutilensuremath {26}&\pgfutilensuremath {24}&\pgfutilensuremath {127}&IMS Alltracks\\%
\pgfutilensuremath {15}&\pgfutilensuremath {9.12}&\pgfutilensuremath {9.12}&\pgfutilensuremath {9.12}&\pgfutilensuremath {9.12}&\pgfutilensuremath {24}&\pgfutilensuremath {24}&\pgfutilensuremath {127}&IMS Alltracks\\%
\pgfutilensuremath {16}&\pgfutilensuremath {9.44}&\pgfutilensuremath {9.44}&\pgfutilensuremath {9.44}&\pgfutilensuremath {9.44}&\pgfutilensuremath {19}&\pgfutilensuremath {22}&\pgfutilensuremath {127}&IMS Alltracks\\%
\pgfutilensuremath {17}&\pgfutilensuremath {9.24}&\pgfutilensuremath {9.24}&\pgfutilensuremath {9.24}&\pgfutilensuremath {9.24}&\pgfutilensuremath {26}&\pgfutilensuremath {32}&\pgfutilensuremath {127}&IMS Alltracks\\%
\pgfutilensuremath {18}&\pgfutilensuremath {8.39}&\pgfutilensuremath {8.39}&\pgfutilensuremath {8.39}&\pgfutilensuremath {8.39}&\pgfutilensuremath {36}&\pgfutilensuremath {33}&\pgfutilensuremath {127}&NYCC\\%
\pgfutilensuremath {19}&\pgfutilensuremath {8.6}&\pgfutilensuremath {8.6}&\pgfutilensuremath {8.6}&\pgfutilensuremath {8.61}&\pgfutilensuremath {32}&\pgfutilensuremath {35}&\pgfutilensuremath {127}&NYCC\\%
\pgfutilensuremath {20}&\pgfutilensuremath {8.78}&\pgfutilensuremath {8.78}&\pgfutilensuremath {8.78}&\pgfutilensuremath {8.78}&\pgfutilensuremath {27}&\pgfutilensuremath {29}&\pgfutilensuremath {127}&NYCC\\%
\pgfutilensuremath {21}&\pgfutilensuremath {8.96}&\pgfutilensuremath {8.96}&\pgfutilensuremath {8.96}&\pgfutilensuremath {8.96}&\pgfutilensuremath {28}&\pgfutilensuremath {28}&\pgfutilensuremath {127}&NYCC\\%
\pgfutilensuremath {22}&\pgfutilensuremath {8.63}&\pgfutilensuremath {8.63}&\pgfutilensuremath {8.63}&\pgfutilensuremath {8.63}&\pgfutilensuremath {35}&\pgfutilensuremath {36}&\pgfutilensuremath {127}&NYCC\\%
\pgfutilensuremath {23}&\pgfutilensuremath {8.91}&\pgfutilensuremath {8.91}&\pgfutilensuremath {8.91}&\pgfutilensuremath {8.92}&\pgfutilensuremath {27}&\pgfutilensuremath {29}&\pgfutilensuremath {127}&NYCC\\%
\pgfutilensuremath {24}&\pgfutilensuremath {9.26}&\pgfutilensuremath {9.26}&\pgfutilensuremath {9.26}&\pgfutilensuremath {9.27}&\pgfutilensuremath {21}&\pgfutilensuremath {16}&\pgfutilensuremath {127}&WLTC\\%
\pgfutilensuremath {25}&\pgfutilensuremath {8.93}&\pgfutilensuremath {8.93}&\pgfutilensuremath {8.93}&\pgfutilensuremath {8.93}&\pgfutilensuremath {34}&\pgfutilensuremath {22}&\pgfutilensuremath {127}&WLTC\\%
\pgfutilensuremath {26}&\pgfutilensuremath {9.45}&\pgfutilensuremath {9.45}&\pgfutilensuremath {9.45}&\pgfutilensuremath {9.45}&\pgfutilensuremath {17}&\pgfutilensuremath {22}&\pgfutilensuremath {127}&WLTC\\%
\pgfutilensuremath {27}&\pgfutilensuremath {9.41}&\pgfutilensuremath {9.41}&\pgfutilensuremath {9.41}&\pgfutilensuremath {9.42}&\pgfutilensuremath {22}&\pgfutilensuremath {18}&\pgfutilensuremath {127}&WLTC\\%
\pgfutilensuremath {28}&\pgfutilensuremath {9.39}&\pgfutilensuremath {9.39}&\pgfutilensuremath {9.39}&\pgfutilensuremath {9.39}&\pgfutilensuremath {27}&\pgfutilensuremath {27}&\pgfutilensuremath {127}&WLTC\\%
\pgfutilensuremath {29}&\pgfutilensuremath {9.35}&\pgfutilensuremath {9.35}&\pgfutilensuremath {9.35}&\pgfutilensuremath {9.36}&\pgfutilensuremath {23}&\pgfutilensuremath {20}&\pgfutilensuremath {127}&WLTC\\%
\pgfutilensuremath {30}&\pgfutilensuremath {9.87}&\pgfutilensuremath {9.87}&\pgfutilensuremath {9.87}&\pgfutilensuremath {9.87}&\pgfutilensuremath {60}&\pgfutilensuremath {62}&\pgfutilensuremath {511}&Artemis Motorway\\%
\pgfutilensuremath {31}&\pgfutilensuremath {9.94}&\pgfutilensuremath {9.94}&\pgfutilensuremath {9.94}&\pgfutilensuremath {9.94}&\pgfutilensuremath {59}&\pgfutilensuremath {58}&\pgfutilensuremath {511}&Artemis Motorway\\%
\pgfutilensuremath {32}&\pgfutilensuremath {10.02}&\pgfutilensuremath {10.02}&\pgfutilensuremath {10.02}&\pgfutilensuremath {10.02}&\pgfutilensuremath {86}&\pgfutilensuremath {84}&\pgfutilensuremath {511}&Artemis Motorway\\%
\pgfutilensuremath {33}&\pgfutilensuremath {10.08}&\pgfutilensuremath {10.08}&\pgfutilensuremath {10.08}&\pgfutilensuremath {10.09}&\pgfutilensuremath {59}&\pgfutilensuremath {59}&\pgfutilensuremath {511}&Artemis Motorway\\%
\pgfutilensuremath {34}&\pgfutilensuremath {9.97}&\pgfutilensuremath {9.97}&\pgfutilensuremath {9.97}&\pgfutilensuremath {9.97}&\pgfutilensuremath {41}&\pgfutilensuremath {47}&\pgfutilensuremath {511}&Artemis Motorway\\%
\pgfutilensuremath {35}&\pgfutilensuremath {9.05}&\pgfutilensuremath {9.05}&\pgfutilensuremath {9.05}&\pgfutilensuremath {9.05}&\pgfutilensuremath {55}&\pgfutilensuremath {59}&\pgfutilensuremath {511}&FTP75\\%
\pgfutilensuremath {36}&\pgfutilensuremath {9.03}&\pgfutilensuremath {9.03}&\pgfutilensuremath {9.03}&\pgfutilensuremath {9.03}&\pgfutilensuremath {39}&\pgfutilensuremath {41}&\pgfutilensuremath {511}&FTP75\\%
\pgfutilensuremath {37}&\pgfutilensuremath {9.23}&\pgfutilensuremath {9.23}&\pgfutilensuremath {9.23}&\pgfutilensuremath {9.24}&\pgfutilensuremath {53}&\pgfutilensuremath {57}&\pgfutilensuremath {511}&FTP75\\%
\pgfutilensuremath {38}&\pgfutilensuremath {9.13}&\pgfutilensuremath {9.13}&\pgfutilensuremath {9.13}&\pgfutilensuremath {9.13}&\pgfutilensuremath {37}&\pgfutilensuremath {39}&\pgfutilensuremath {511}&FTP75\\%
\pgfutilensuremath {39}&\pgfutilensuremath {8.81}&\pgfutilensuremath {8.81}&\pgfutilensuremath {8.81}&\pgfutilensuremath {8.81}&\pgfutilensuremath {63}&\pgfutilensuremath {70}&\pgfutilensuremath {511}&FTP75\\%
\pgfutilensuremath {40}&\pgfutilensuremath {9.34}&\pgfutilensuremath {9.34}&\pgfutilensuremath {9.34}&\pgfutilensuremath {9.34}&\pgfutilensuremath {45}&\pgfutilensuremath {50}&\pgfutilensuremath {511}&FTP75\\%
\pgfutilensuremath {41}&\pgfutilensuremath {9.48}&\pgfutilensuremath {9.48}&\pgfutilensuremath {9.48}&\pgfutilensuremath {9.48}&\pgfutilensuremath {34}&\pgfutilensuremath {35}&\pgfutilensuremath {511}&IMS Alltracks\\%
\pgfutilensuremath {42}&\pgfutilensuremath {9.02}&\pgfutilensuremath {9.02}&\pgfutilensuremath {9.02}&\pgfutilensuremath {9.02}&\pgfutilensuremath {50}&\pgfutilensuremath {51}&\pgfutilensuremath {511}&IMS Alltracks\\%
\pgfutilensuremath {43}&\pgfutilensuremath {9.31}&\pgfutilensuremath {9.31}&\pgfutilensuremath {9.31}&\pgfutilensuremath {9.31}&\pgfutilensuremath {38}&\pgfutilensuremath {43}&\pgfutilensuremath {511}&IMS Alltracks\\%
\pgfutilensuremath {44}&\pgfutilensuremath {9.4}&\pgfutilensuremath {9.4}&\pgfutilensuremath {9.4}&\pgfutilensuremath {9.41}&\pgfutilensuremath {34}&\pgfutilensuremath {33}&\pgfutilensuremath {511}&IMS Alltracks\\%
\pgfutilensuremath {45}&\pgfutilensuremath {9.37}&\pgfutilensuremath {9.37}&\pgfutilensuremath {9.37}&\pgfutilensuremath {9.38}&\pgfutilensuremath {41}&\pgfutilensuremath {41}&\pgfutilensuremath {511}&IMS Alltracks\\%
\pgfutilensuremath {46}&\pgfutilensuremath {8.35}&\pgfutilensuremath {8.35}&\pgfutilensuremath {8.35}&\pgfutilensuremath {8.35}&\pgfutilensuremath {86}&\pgfutilensuremath {94}&\pgfutilensuremath {511}&NYCC\\%
\pgfutilensuremath {47}&\pgfutilensuremath {8.59}&\pgfutilensuremath {8.59}&\pgfutilensuremath {8.59}&\pgfutilensuremath {8.59}&\pgfutilensuremath {68}&\pgfutilensuremath {61}&\pgfutilensuremath {511}&NYCC\\%
\pgfutilensuremath {48}&\pgfutilensuremath {8.86}&\pgfutilensuremath {8.86}&\pgfutilensuremath {8.86}&\pgfutilensuremath {8.86}&\pgfutilensuremath {85}&\pgfutilensuremath {83}&\pgfutilensuremath {511}&NYCC\\%
\pgfutilensuremath {49}&\pgfutilensuremath {8.63}&\pgfutilensuremath {8.63}&\pgfutilensuremath {8.63}&\pgfutilensuremath {8.63}&\pgfutilensuremath {65}&\pgfutilensuremath {61}&\pgfutilensuremath {511}&NYCC\\%
\pgfutilensuremath {50}&\pgfutilensuremath {8.93}&\pgfutilensuremath {8.93}&\pgfutilensuremath {8.93}&\pgfutilensuremath {8.93}&\pgfutilensuremath {54}&\pgfutilensuremath {53}&\pgfutilensuremath {511}&NYCC\\%
\pgfutilensuremath {51}&\pgfutilensuremath {8.56}&\pgfutilensuremath {8.56}&\pgfutilensuremath {8.56}&\pgfutilensuremath {8.56}&\pgfutilensuremath {66}&\pgfutilensuremath {66}&\pgfutilensuremath {511}&NYCC\\%
\pgfutilensuremath {52}&\pgfutilensuremath {9.27}&\pgfutilensuremath {9.27}&\pgfutilensuremath {9.27}&\pgfutilensuremath {9.27}&\pgfutilensuremath {56}&\pgfutilensuremath {61}&\pgfutilensuremath {511}&WLTC\\%
\pgfutilensuremath {53}&\pgfutilensuremath {9.41}&\pgfutilensuremath {9.41}&\pgfutilensuremath {9.41}&\pgfutilensuremath {9.42}&\pgfutilensuremath {58}&\pgfutilensuremath {60}&\pgfutilensuremath {511}&WLTC\\%
\pgfutilensuremath {54}&\pgfutilensuremath {9.55}&\pgfutilensuremath {9.55}&\pgfutilensuremath {9.55}&\pgfutilensuremath {9.56}&\pgfutilensuremath {47}&\pgfutilensuremath {44}&\pgfutilensuremath {511}&WLTC\\%
\pgfutilensuremath {55}&\pgfutilensuremath {9.43}&\pgfutilensuremath {9.43}&\pgfutilensuremath {9.43}&\pgfutilensuremath {9.43}&\pgfutilensuremath {51}&\pgfutilensuremath {61}&\pgfutilensuremath {511}&WLTC\\%
\pgfutilensuremath {56}&\pgfutilensuremath {9.51}&\pgfutilensuremath {9.51}&\pgfutilensuremath {9.51}&\pgfutilensuremath {9.51}&\pgfutilensuremath {44}&\pgfutilensuremath {61}&\pgfutilensuremath {511}&WLTC\\%
\pgfutilensuremath {57}&\pgfutilensuremath {9.58}&\pgfutilensuremath {9.58}&\pgfutilensuremath {9.58}&\pgfutilensuremath {9.58}&\pgfutilensuremath {38}&\pgfutilensuremath {41}&\pgfutilensuremath {511}&WLTC\\%
\pgfutilensuremath {58}&\pgfutilensuremath {9.86}&\pgfutilensuremath {9.86}&\pgfutilensuremath {9.86}&\pgfutilensuremath {9.86}&\pgfutilensuremath {102}&\pgfutilensuremath {105}&\pgfutilensuremath {2{,}047}&Artemis Motorway\\%
\pgfutilensuremath {59}&\pgfutilensuremath {9.77}&\pgfutilensuremath {9.77}&\pgfutilensuremath {9.77}&\pgfutilensuremath {9.77}&\pgfutilensuremath {133}&\pgfutilensuremath {128}&\pgfutilensuremath {2{,}047}&Artemis Motorway\\%
\pgfutilensuremath {60}&\pgfutilensuremath {9.96}&\pgfutilensuremath {9.96}&\pgfutilensuremath {9.96}&\pgfutilensuremath {9.96}&\pgfutilensuremath {108}&\pgfutilensuremath {113}&\pgfutilensuremath {2{,}047}&Artemis Motorway\\%
\pgfutilensuremath {61}&\pgfutilensuremath {9.46}&\pgfutilensuremath {9.46}&\pgfutilensuremath {9.46}&\pgfutilensuremath {9.46}&\pgfutilensuremath {18}&\pgfutilensuremath {22}&\pgfutilensuremath {127}&IMS Alltracks\\%
\pgfutilensuremath {62}&\pgfutilensuremath {9.75}&\pgfutilensuremath {9.75}&\pgfutilensuremath {9.75}&\pgfutilensuremath {9.75}&\pgfutilensuremath {58}&\pgfutilensuremath {57}&\pgfutilensuremath {511}&Artemis Motorway\\%
\pgfutilensuremath {63}&\pgfutilensuremath {9.19}&\pgfutilensuremath {9.19}&\pgfutilensuremath {9.19}&\pgfutilensuremath {9.19}&\pgfutilensuremath {45}&\pgfutilensuremath {50}&\pgfutilensuremath {511}&IMS Alltracks\\%
\pgfutilensuremath {64}&\pgfutilensuremath {9.93}&\pgfutilensuremath {9.93}&\pgfutilensuremath {9.93}&\pgfutilensuremath {9.93}&\pgfutilensuremath {130}&\pgfutilensuremath {130}&\pgfutilensuremath {2{,}047}&Artemis Motorway\\%
\pgfutilensuremath {65}&\pgfutilensuremath {9.05}&\pgfutilensuremath {9.05}&\pgfutilensuremath {9.05}&\pgfutilensuremath {9.06}&\pgfutilensuremath {169}&\pgfutilensuremath {161}&\pgfutilensuremath {2{,}047}&NYCC\\%
\pgfutilensuremath {66}&\pgfutilensuremath {9.67}&\pgfutilensuremath {9.67}&\pgfutilensuremath {9.67}&\pgfutilensuremath {9.67}&\pgfutilensuremath {8}&\pgfutilensuremath {8}&\pgfutilensuremath {31}&Artemis Motorway\\%
\pgfutilensuremath {67}&\pgfutilensuremath {10.21}&\pgfutilensuremath {10.21}&\pgfutilensuremath {10.21}&\pgfutilensuremath {10.21}&\pgfutilensuremath {132}&\pgfutilensuremath {133}&\pgfutilensuremath {2{,}047}&Artemis Motorway\\%
\pgfutilensuremath {68}&\pgfutilensuremath {10.08}&\pgfutilensuremath {10.08}&\pgfutilensuremath {10.08}&\pgfutilensuremath {10.08}&\pgfutilensuremath {126}&\pgfutilensuremath {130}&\pgfutilensuremath {2{,}047}&Artemis Motorway\\%
\pgfutilensuremath {69}&\pgfutilensuremath {9.32}&\pgfutilensuremath {9.32}&\pgfutilensuremath {9.32}&\pgfutilensuremath {9.32}&\pgfutilensuremath {112}&\pgfutilensuremath {110}&\pgfutilensuremath {2{,}047}&FTP75\\%
\pgfutilensuremath {70}&\pgfutilensuremath {8.9}&\pgfutilensuremath {8.9}&\pgfutilensuremath {8.9}&\pgfutilensuremath {8.9}&\pgfutilensuremath {112}&\pgfutilensuremath {127}&\pgfutilensuremath {2{,}047}&FTP75\\%
\pgfutilensuremath {71}&\pgfutilensuremath {9.2}&\pgfutilensuremath {9.2}&\pgfutilensuremath {9.2}&\pgfutilensuremath {9.2}&\pgfutilensuremath {97}&\pgfutilensuremath {98}&\pgfutilensuremath {2{,}047}&FTP75\\%
\pgfutilensuremath {72}&\pgfutilensuremath {9.27}&\pgfutilensuremath {9.27}&\pgfutilensuremath {9.27}&\pgfutilensuremath {9.27}&\pgfutilensuremath {133}&\pgfutilensuremath {146}&\pgfutilensuremath {2{,}047}&FTP75\\%
\pgfutilensuremath {73}&\pgfutilensuremath {9.37}&\pgfutilensuremath {9.37}&\pgfutilensuremath {9.37}&\pgfutilensuremath {9.37}&\pgfutilensuremath {119}&\pgfutilensuremath {120}&\pgfutilensuremath {2{,}047}&FTP75\\%
\pgfutilensuremath {74}&\pgfutilensuremath {9.29}&\pgfutilensuremath {9.29}&\pgfutilensuremath {9.29}&\pgfutilensuremath {9.29}&\pgfutilensuremath {88}&\pgfutilensuremath {81}&\pgfutilensuremath {2{,}047}&FTP75\\%
\pgfutilensuremath {75}&\pgfutilensuremath {9.26}&\pgfutilensuremath {9.26}&\pgfutilensuremath {9.26}&\pgfutilensuremath {9.26}&\pgfutilensuremath {85}&\pgfutilensuremath {88}&\pgfutilensuremath {2{,}047}&IMS Alltracks\\%
\pgfutilensuremath {76}&\pgfutilensuremath {9.3}&\pgfutilensuremath {9.3}&\pgfutilensuremath {9.3}&\pgfutilensuremath {9.31}&\pgfutilensuremath {83}&\pgfutilensuremath {77}&\pgfutilensuremath {2{,}047}&IMS Alltracks\\%
\pgfutilensuremath {77}&\pgfutilensuremath {9.54}&\pgfutilensuremath {9.54}&\pgfutilensuremath {9.54}&\pgfutilensuremath {9.54}&\pgfutilensuremath {86}&\pgfutilensuremath {91}&\pgfutilensuremath {2{,}047}&IMS Alltracks\\%
\pgfutilensuremath {78}&\pgfutilensuremath {9.4}&\pgfutilensuremath {9.4}&\pgfutilensuremath {9.4}&\pgfutilensuremath {9.4}&\pgfutilensuremath {67}&\pgfutilensuremath {60}&\pgfutilensuremath {2{,}047}&IMS Alltracks\\%
\pgfutilensuremath {79}&\pgfutilensuremath {9.7}&\pgfutilensuremath {9.7}&\pgfutilensuremath {9.7}&\pgfutilensuremath {9.7}&\pgfutilensuremath {69}&\pgfutilensuremath {67}&\pgfutilensuremath {2{,}047}&IMS Alltracks\\%
\pgfutilensuremath {80}&\pgfutilensuremath {9.59}&\pgfutilensuremath {9.59}&\pgfutilensuremath {9.59}&\pgfutilensuremath {9.59}&\pgfutilensuremath {95}&\pgfutilensuremath {95}&\pgfutilensuremath {2{,}047}&IMS Alltracks\\%
\pgfutilensuremath {81}&\pgfutilensuremath {8.83}&\pgfutilensuremath {8.83}&\pgfutilensuremath {8.83}&\pgfutilensuremath {8.84}&\pgfutilensuremath {162}&\pgfutilensuremath {170}&\pgfutilensuremath {2{,}047}&NYCC\\%
\pgfutilensuremath {82}&\pgfutilensuremath {8.91}&\pgfutilensuremath {8.91}&\pgfutilensuremath {8.91}&\pgfutilensuremath {8.91}&\pgfutilensuremath {147}&\pgfutilensuremath {142}&\pgfutilensuremath {2{,}047}&NYCC\\%
\pgfutilensuremath {83}&\pgfutilensuremath {9.03}&\pgfutilensuremath {9.03}&\pgfutilensuremath {9.03}&\pgfutilensuremath {9.03}&\pgfutilensuremath {165}&\pgfutilensuremath {148}&\pgfutilensuremath {2{,}047}&NYCC\\%
\pgfutilensuremath {84}&\pgfutilensuremath {8.99}&\pgfutilensuremath {8.99}&\pgfutilensuremath {8.99}&\pgfutilensuremath {8.99}&\pgfutilensuremath {140}&\pgfutilensuremath {138}&\pgfutilensuremath {2{,}047}&NYCC\\%
\pgfutilensuremath {85}&\pgfutilensuremath {9.06}&\pgfutilensuremath {9.06}&\pgfutilensuremath {9.06}&\pgfutilensuremath {9.06}&\pgfutilensuremath {176}&\pgfutilensuremath {164}&\pgfutilensuremath {2{,}047}&NYCC\\%
\pgfutilensuremath {86}&\pgfutilensuremath {9.46}&\pgfutilensuremath {9.46}&\pgfutilensuremath {9.46}&\pgfutilensuremath {9.46}&\pgfutilensuremath {95}&\pgfutilensuremath {98}&\pgfutilensuremath {2{,}047}&WLTC\\%
\pgfutilensuremath {87}&\pgfutilensuremath {9.58}&\pgfutilensuremath {9.58}&\pgfutilensuremath {9.58}&\pgfutilensuremath {9.58}&\pgfutilensuremath {108}&\pgfutilensuremath {108}&\pgfutilensuremath {2{,}047}&WLTC\\%
\pgfutilensuremath {88}&\pgfutilensuremath {9.66}&\pgfutilensuremath {9.66}&\pgfutilensuremath {9.66}&\pgfutilensuremath {9.66}&\pgfutilensuremath {106}&\pgfutilensuremath {104}&\pgfutilensuremath {2{,}047}&WLTC\\%
\pgfutilensuremath {89}&\pgfutilensuremath {9.51}&\pgfutilensuremath {9.51}&\pgfutilensuremath {9.51}&\pgfutilensuremath {9.51}&\pgfutilensuremath {97}&\pgfutilensuremath {99}&\pgfutilensuremath {2{,}047}&WLTC\\%
\pgfutilensuremath {90}&\pgfutilensuremath {9.63}&\pgfutilensuremath {9.63}&\pgfutilensuremath {9.63}&\pgfutilensuremath {9.63}&\pgfutilensuremath {111}&\pgfutilensuremath {115}&\pgfutilensuremath {2{,}047}&WLTC\\%
\pgfutilensuremath {91}&\pgfutilensuremath {9.72}&\pgfutilensuremath {9.72}&\pgfutilensuremath {9.72}&\pgfutilensuremath {9.72}&\pgfutilensuremath {157}&\pgfutilensuremath {152}&\pgfutilensuremath {2{,}047}&WLTC\\%
\pgfutilensuremath {92}&\pgfutilensuremath {9.58}&\pgfutilensuremath {9.58}&\pgfutilensuremath {9.58}&\pgfutilensuremath {9.59}&\pgfutilensuremath {8}&\pgfutilensuremath {8}&\pgfutilensuremath {31}&Artemis Motorway\\%
\pgfutilensuremath {93}&\pgfutilensuremath {9.1}&\pgfutilensuremath {9.1}&\pgfutilensuremath {9.1}&\pgfutilensuremath {9.1}&\pgfutilensuremath {12}&\pgfutilensuremath {14}&\pgfutilensuremath {31}&Artemis Motorway\\%
\pgfutilensuremath {94}&\pgfutilensuremath {8.98}&\pgfutilensuremath {8.98}&\pgfutilensuremath {8.98}&\pgfutilensuremath {8.99}&\pgfutilensuremath {15}&\pgfutilensuremath {8}&\pgfutilensuremath {31}&Artemis Motorway\\%
\pgfutilensuremath {95}&\pgfutilensuremath {9.06}&\pgfutilensuremath {9.06}&\pgfutilensuremath {9.06}&\pgfutilensuremath {9.06}&\pgfutilensuremath {12}&\pgfutilensuremath {11}&\pgfutilensuremath {31}&Artemis Motorway\\%
\pgfutilensuremath {96}&\pgfutilensuremath {9.48}&\pgfutilensuremath {9.48}&\pgfutilensuremath {9.48}&\pgfutilensuremath {9.48}&\pgfutilensuremath {8}&\pgfutilensuremath {7}&\pgfutilensuremath {31}&Artemis Motorway\\%
\pgfutilensuremath {97}&\pgfutilensuremath {8.09}&\pgfutilensuremath {8.09}&\pgfutilensuremath {8.09}&\pgfutilensuremath {8.09}&\pgfutilensuremath {16}&\pgfutilensuremath {17}&\pgfutilensuremath {31}&FTP75\\%
\pgfutilensuremath {98}&\pgfutilensuremath {8.98}&\pgfutilensuremath {8.98}&\pgfutilensuremath {8.98}&\pgfutilensuremath {8.98}&\pgfutilensuremath {10}&\pgfutilensuremath {11}&\pgfutilensuremath {31}&FTP75\\%
\pgfutilensuremath {99}&\pgfutilensuremath {9.16}&\pgfutilensuremath {9.16}&\pgfutilensuremath {9.16}&\pgfutilensuremath {9.19}&\pgfutilensuremath {9}&\pgfutilensuremath {9}&\pgfutilensuremath {31}&FTP75\\%
\pgfutilensuremath {100}&\pgfutilensuremath {8.84}&\pgfutilensuremath {8.84}&\pgfutilensuremath {8.84}&\pgfutilensuremath {8.84}&\pgfutilensuremath {9}&\pgfutilensuremath {8}&\pgfutilensuremath {31}&FTP75\\%
\pgfutilensuremath {101}&\pgfutilensuremath {8.82}&\pgfutilensuremath {8.82}&\pgfutilensuremath {8.82}&\pgfutilensuremath {8.84}&\pgfutilensuremath {9}&\pgfutilensuremath {7}&\pgfutilensuremath {31}&FTP75\\%
\pgfutilensuremath {102}&\pgfutilensuremath {9.04}&\pgfutilensuremath {9.04}&\pgfutilensuremath {9.04}&\pgfutilensuremath {9.04}&\pgfutilensuremath {8}&\pgfutilensuremath {9}&\pgfutilensuremath {31}&FTP75\\%
\pgfutilensuremath {103}&\pgfutilensuremath {8.65}&\pgfutilensuremath {8.65}&\pgfutilensuremath {8.65}&\pgfutilensuremath {8.65}&\pgfutilensuremath {15}&\pgfutilensuremath {13}&\pgfutilensuremath {31}&IMS Alltracks\\%
\pgfutilensuremath {104}&\pgfutilensuremath {9.24}&\pgfutilensuremath {9.24}&\pgfutilensuremath {9.24}&\pgfutilensuremath {9.24}&\pgfutilensuremath {10}&\pgfutilensuremath {10}&\pgfutilensuremath {31}&IMS Alltracks\\%
\pgfutilensuremath {105}&\pgfutilensuremath {9.17}&\pgfutilensuremath {9.17}&\pgfutilensuremath {9.17}&\pgfutilensuremath {9.17}&\pgfutilensuremath {8}&\pgfutilensuremath {9}&\pgfutilensuremath {31}&IMS Alltracks\\%
\pgfutilensuremath {106}&\pgfutilensuremath {8.93}&\pgfutilensuremath {8.93}&\pgfutilensuremath {8.93}&\pgfutilensuremath {8.93}&\pgfutilensuremath {12}&\pgfutilensuremath {9}&\pgfutilensuremath {31}&IMS Alltracks\\%
\pgfutilensuremath {107}&\pgfutilensuremath {8.77}&\pgfutilensuremath {8.77}&\pgfutilensuremath {8.77}&\pgfutilensuremath {8.77}&\pgfutilensuremath {13}&\pgfutilensuremath {13}&\pgfutilensuremath {31}&IMS Alltracks\\%
\pgfutilensuremath {108}&\pgfutilensuremath {9.13}&\pgfutilensuremath {9.13}&\pgfutilensuremath {9.13}&\pgfutilensuremath {9.13}&\pgfutilensuremath {10}&\pgfutilensuremath {14}&\pgfutilensuremath {31}&IMS Alltracks\\%
\pgfutilensuremath {109}&\pgfutilensuremath {8.56}&\pgfutilensuremath {8.56}&\pgfutilensuremath {8.56}&\pgfutilensuremath {8.56}&\pgfutilensuremath {8}&\pgfutilensuremath {8}&\pgfutilensuremath {31}&NYCC\\%
\pgfutilensuremath {110}&\pgfutilensuremath {8.54}&\pgfutilensuremath {8.54}&\pgfutilensuremath {8.54}&\pgfutilensuremath {8.54}&\pgfutilensuremath {9}&\pgfutilensuremath {10}&\pgfutilensuremath {31}&NYCC\\%
\pgfutilensuremath {111}&\pgfutilensuremath {8.09}&\pgfutilensuremath {8.09}&\pgfutilensuremath {8.09}&\pgfutilensuremath {8.09}&\pgfutilensuremath {9}&\pgfutilensuremath {10}&\pgfutilensuremath {31}&NYCC\\%
\pgfutilensuremath {112}&\pgfutilensuremath {7.89}&\pgfutilensuremath {7.89}&\pgfutilensuremath {7.89}&\pgfutilensuremath {7.89}&\pgfutilensuremath {12}&\pgfutilensuremath {13}&\pgfutilensuremath {31}&NYCC\\%
\pgfutilensuremath {113}&\pgfutilensuremath {8.54}&\pgfutilensuremath {8.54}&\pgfutilensuremath {8.54}&\pgfutilensuremath {8.54}&\pgfutilensuremath {9}&\pgfutilensuremath {13}&\pgfutilensuremath {31}&NYCC\\%
\pgfutilensuremath {114}&\pgfutilensuremath {8.75}&\pgfutilensuremath {8.75}&\pgfutilensuremath {8.75}&\pgfutilensuremath {8.75}&\pgfutilensuremath {8}&\pgfutilensuremath {8}&\pgfutilensuremath {31}&NYCC\\%
\pgfutilensuremath {115}&\pgfutilensuremath {9.08}&\pgfutilensuremath {9.08}&\pgfutilensuremath {9.08}&\pgfutilensuremath {9.08}&\pgfutilensuremath {8}&\pgfutilensuremath {9}&\pgfutilensuremath {31}&WLTC\\%
\pgfutilensuremath {116}&\pgfutilensuremath {8.88}&\pgfutilensuremath {8.88}&\pgfutilensuremath {8.88}&\pgfutilensuremath {8.89}&\pgfutilensuremath {7}&\pgfutilensuremath {7}&\pgfutilensuremath {31}&WLTC\\%
\pgfutilensuremath {117}&\pgfutilensuremath {9.1}&\pgfutilensuremath {9.1}&\pgfutilensuremath {9.1}&\pgfutilensuremath {9.1}&\pgfutilensuremath {10}&\pgfutilensuremath {9}&\pgfutilensuremath {31}&WLTC\\%
\pgfutilensuremath {118}&\pgfutilensuremath {9.21}&\pgfutilensuremath {9.21}&\pgfutilensuremath {9.21}&\pgfutilensuremath {9.21}&\pgfutilensuremath {9}&\pgfutilensuremath {7}&\pgfutilensuremath {31}&WLTC\\%
\pgfutilensuremath {119}&\pgfutilensuremath {8.64}&\pgfutilensuremath {8.64}&\pgfutilensuremath {8.64}&\pgfutilensuremath {8.64}&\pgfutilensuremath {12}&\pgfutilensuremath {10}&\pgfutilensuremath {31}&WLTC\\%
\pgfutilensuremath {120}&\pgfutilensuremath {8.64}&\pgfutilensuremath {8.64}&\pgfutilensuremath {8.64}&\pgfutilensuremath {8.64}&\pgfutilensuremath {10}&\pgfutilensuremath {11}&\pgfutilensuremath {31}&WLTC\\%
\label{tab:a1}
\end{longtable}%

\end{scriptsize}